\input amstex
\documentstyle{amsppt}
\magnification 1200
\vcorrection{-9mm}
\NoBlackBoxes
\input epsf

\def\Z{\Bbb Z}
\def\R{\Bbb R}
\def\C{\Bbb C}
\def\P{\Bbb P}
\def\sph{\Bbb S}
\def\RP{\Bbb{RP}}
\def\CP{\Bbb{CP}}
\def\conj{\operatorname{conj}}

\def\hyp{\operatorname{hyp}}
\def\ord{\operatorname{ord}}

\def\card{\operatorname{Card}}

\def\Sing{\operatorname{Sing}}

\topmatter
\title    Topology of maximally writhed real algebraic knots
\endtitle

\author   G.~B.~Mikhalkin and S.~Yu.~Orevkov
\endauthor

\abstract
             Oleg Viro introduced an invariant of rigid isotopy for real
             algebraic knots in $\Bbb{RP}^3$ which can be viewed
             as a first order Vassiliev invariant.
             In this paper we look at real algebraic knots of degree $d$ with the maximal possible value
             of this invariant. We show that for a given $d$ all such knots are topologically isotopic and
             explicitly identify their knot type.
%
%             We prove that if for a knot
%             of degree $d$, this invariant attains the maximal possible value, % $N_d:=(d-1)(d-2)/2$,
%             then the knot is covered by the $(d,d-2)$-torus link under the covering
%             $\sph^3\to\Bbb{RP}^3$.
%             We also construct $[d/2]-1$ pairwise non-isotopic real algebraic knots of any degree $d$
%             with Viro invariant $N_d-2$.
\endabstract

\thanks
Research is supported in part by  the SNSF-grants
159240, 159581 and the NCCR SwissMAP project  (G.M.)
and by RSF grant, project 14-21-00053 dated 11.08.14. (S.O.)
\endthanks

\address        Universit\'e de Gen\'eve, Section de Math\'ematiques,
                Battelle Villa, 1227 Carouge, Suisse
\endaddress
\email          grigory.mikhalkin\@unige.ch
\endemail

\address
                Steklov Mathematical Institute, Gubkina 8, 119991, Moscow, Russia;\smallskip
                IMT, Universit\'e Paul Sabatier,
                118 route de Narbonne, 31062, Toulouse, France;\smallskip
                National Research University Higher School of Economics, Russian Federation,
                Usacheva 6, Moscow, 119048 Russia
\endaddress
%-----
\email          orevkov\@math.ups-tlse.fr
\endemail
%-----

\endtopmatter

\def\refB   {1}
\def\refMO  {2}
\def\refMu  {3}
\def\refV   {4}

\def\sectTors        {1}
\def\sectProofMainTh {2}

\def\lemTors        {1}

\def\lemHypC        {2}
\def\lemHypKC       {3}
\def\lemSmooth      {4}
\def\lemHyperb      {5}
\def\lemHypK        {6}
\def\lemDisjoint    {7}
\def\lemTwoLines    {8}
\def\propU          {2}

\def\eqDefHyp   {1}
\def\eqIndex    {2}

\def\figCusp    {1}
\def\figR       {2}

\document

%A {\it link} in a 3-manifold is a disjoint union of smoothly embedded circles considered up
%to isotopy. A one-component link is called {\it knot}.

\head Introduction
\endhead

A {\it real algebraic curve} in $\P^3$ is a (complex) one-dimensional subvariety $L$ in
$\P^3=\CP^3$ invariant under the involution of complex conjugation
$\conj:\P^3\to\P^3$, $(x_0:x_1:x_2:x_3)\mapsto(\bar x_0:\bar x_1:\bar x_2:\bar x_3)$.
The conj-invariance is equivalent to the fact that $L$ can be defined by a system of
homogeneous polynomial equations with real coefficients.
The degree of $L$ is defined as its homological degree, i.~e. the number $d$ such that
$[L] = d[\P^1]\in H_2(\P^3)\cong\Z$. A curve of degree $d$ intersects a generic complex plane
in $d$ points.

We denote the set of real points of $L$ by $\R L$.
We say that a real curve $L$ is {\it smooth} if it is a smooth complex submanifold of $\P^3$.
In this case, $\R L$ is a smooth real submanifold of $\RP^3$ and if it is non-empty,
%we call it {\it real algebraic link} (or {\it knot} when $\R L$ is connected).
we call it a {\it real algebraic link} or, more specifically,
a {\it real algebraic knot} in the case when $\R L$ is connected.

Two real algebraic links are called {\it rigidly isotopic} if they belong to the same connected
component of the space of smooth real curves of the same degree.
A rigid isotopy classification of real algebraic rational curves in $\P^3$ is obtained in [\refB] up to degree 5
and in [\refMO] up to degree 6. Also we gave in [\refMO] a rigid isotopy classification for genus one knots
and links up to degree 6 (here we speak of the genus of the complex curve $L$ rather than the
minimal genus of a Seifert surface of $\R L$).

In all the above-mentioned cases, a rigid isotopy class is
completely determined by the usual (topological) isotopy class, the complex orientation (for genus one links),
and the invariant of rigid isotopy $w$ introduced by Viro [\refV] (called in [\refV] {\it encomplexed writhe\/}).
This invariant is defined as the sum of signs of crossings of a generic projection but the
crossings with non-real branches are also counted with appropriate signs; see details in [\refV]
(the definition of $w$ is also reproduced in [\refMO]).

Let $T(p,q)=\{(z,w)\mid z^p=w^q\}\cap \sph^3$, $p\ge q\ge 0$,
be the $(p,q)$-torus link in the 3-sphere $\sph^3\subset\C^2$.
If $p\equiv q\mod2$, we define the {\it projective torus link} %(or {\it link\/})
$T_{proj}(p,q) = T(p,q)/(-1) \subset \sph^3/(-1) = \RP^3$.

\smallskip
Let $N_d = (d-1)(d-2)/2$. By the genus formula, this is the maximal possible value of $w$ for
irreducible curves of degree $d$ which can be attained on rational curves only.
So, if a real algebraic curve $K$ in $\P^3$ is smooth, irreducible, and $|w(K)|=N_d$ where
$d=\deg K$ (and hence the genus of $K$ is zero), then we call it {\it maximally writhed} or {\it $MW$-curve}.
The main result is the following.

\proclaim{ Theorem 1 }
Let $K$ be an $MW$-curve of degree $d\ge 3$, and $w(K)=N_d$.
Then $\R K$ is isotopic to $T_{proj}(d,d-2)$.
\endproclaim

%We say that a real algebraic knot is {\it maximally writhed} or {\it $MW$-knot}
%if it is as in Theorem 1.

\proclaim{ Corollary 1 }
A plane projection of an $MW$-curve from any generic real point
%
%Let $K$ be an $MW$-knot of degree $d$.
%Then any generic plane projection of $\R K$ 
has $N_d$ or $N_d-1$
real double points with real local branches.
\endproclaim

\demo{ Proof } Follows from Murasugi's result [\refMu; Proposition 7.5] which states that
any projection of a torus link $T(p,q)$, $1\le q\le p$, has at least $p(q-1)$ crossings. \qed
\enddemo

\noindent
In Proposition \propU\ (see the end of the paper) we give a precision and
a self-contained (i.e., not using [\refMu]) proof to Corollary 1.

\smallskip
\noindent
{\bf Conjecture 1. } In Theorem 1, $K$ is rigidly isotopic to $T_{proj}(d,d-2)$.
%Any $MW$-knot of degree $d$ is rigidly isotopic to $T_{proj}(d,d-2)$.

\smallskip
\noindent
{\bf Conjecture 2. } If an algebraic knot $\R K$ of degree $d$ in $\RP^3$ is isotopic to $T_{proj}(d,d-2)$,
then $w(K)=N_d$.

\smallskip
In a forthcoming paper we are going to give a proof of Conjecture 2 as well as a generalization
of Theorem~1 for links of arbitrary genus.

\smallskip
The following differential geometric property of maximally writhed algebraic knots was
communicated to us by Oleg Viro.

\proclaim{ Proposition 1 }
Let $K$ be as in Theorem 1. Then the torsion of $\R K$ is everywhere positive.
\endproclaim

%%%%%%%%%%%%%%%%%%%%%%%%%%%%%%%%%%%%%%%%%%%%%%%%%%%%%%%%%%%%%%%%%%%%%%%%%%%%%%%%%%%%%%%%%%%%%%%%%%%%%
%%%%%%%%%%%%%%%%%%%%%%%%%%%%%%%%%%%%%%%%%%%%%%%%%%%%%%%%%%%%%%%%%%%%%%%%%%%%%%%%%%%%%%%%%%%%%%%%%%%%%
%%%%%%%%%%%%%%%%%%%%%%%%%%%%%%%%%%%%%%%%%%%%%%%%%%%%%%%%%%%%%%%%%%%%%%%%%%%%%%%%%%%%%%%%%%%%%%%%%%%%%

%%%%%%%%%%%%%%%%%%%%%%%%%%%%%%%%%%%%%%%%%%%%%%%%%%%%%%%%%%%%%%%%%%%%%%%%%%%%%%%%%%%%%%%%%%%%%%%%%%%%%
%%%%%%%%%%%%%%%%%%%%%%%%%%%%%%%%%%%%%%%%%%%%%%%%%%%%%%%%%%%%%%%%%%%%%%%%%%%%%%%%%%%%%%%%%%%%%%%%%%%%%
%%%%%%%%%%%%%%%%%%%%%%%%%%%%%%%%%%%%%%%%%%%%%%%%%%%%%%%%%%%%%%%%%%%%%%%%%%%%%%%%%%%%%%%%%%%%%%%%%%%%%

\head\sectTors. $MW$-curves have everywhere positive torsion (proof of Proposition 1)
\endhead

Recall that the sign of the (differential geometric) torsion of a curve
 $t\mapsto r(t)\in\R^3$, $t\in\R$, coincides with the sign of $\det(r',r'',r''')$
and it does not depend on the parametrization if $r'\neq 0$.
The sign of the torsion of a curve in $\RP^3$ does not depend
on a choice of positively oriented affine chart.

\proclaim{ Lemma \lemTors } Let $K$ be a real algebraic knot in $\P^3$ of genus $0$
which is not contained in any plane.
If the torsion of $\R K$ vanishes at a point $p$, then there exists an
arbitrarily small deformation of $K$ (in the class of real algebraic knots)
which has points with negative torsion.
\endproclaim

\demo{ Proof }
We can always choose affine coordinates $(x_1,x_2,x_3)$ centered at $p$ such that
a parametrization $t\mapsto(x_1(t),x_2(t),x_3(t))$ of $K$ at $p$ satisfies the condition
$\ord_t x_1<\ord_t x_2<\ord_t x_3$.
If $\ord_t x_k>k$ for $k=1,2,$ or $3$, then $x_k$ can be perturbed so that
$\ord_t x_k = k$ and the $k$-th derivative $x_k^{(k)}(0)$ has any sign we want.
Indeed, let $(y_0:y_1:y_2:y_3)$ be homogeneous coordinates such that $x_i=y_i/y_0$, $i=1,2,3$.
Then the parametrization can be chosen so that $x_i(t)=y_i(t)/y_0(t)$, $i=1,2,3$, where
$y_0(t),\dots,y_3(t)$ are real polynomials of degree $d=\deg K$, and $y_0(0)>0$.
Then the desired perturbation of $x_k(t)$ is just $(c_k t^k + y_k(t))/y_0(t)$ where $0<|c_k|\ll1$
and $c_k$ has the prescribed sign.
\qed\enddemo

\demo{ Proof of Proposition 1 }
By Lemma \lemTors, it is enough to show that $\R K$ does not have points with negative torsion.
Suppose it does.
Then, in an appropriate affine chart, %centered at a generic point,
$\R K$ admits a parametrization
of the form $t\mapsto(t,t^2+O(t^3),-t^3+O(t^4))$.
This means that in a sufficiently small neighbourhood of the origin, the curve is approximated by
a negatively twisted rational cubic curve. %Such a curve has a projection with a negative crossing
Hence there is a projection with a negative crossing
(see [\refV; Section 1.4]). Since $w(K)$ is the sum of the signs of all real crossings and the number
of them is at most $N_d$, a single negative crossing makes impossible to
attain the equality $w(K)=N_d$.
\qed\enddemo

%%%%%%%%%%%%%%%%%%%%%%%%%%%%%%%%%%%%%%%%%%%%%%%%%%%%%%%%%%%%%%%%%%%%%%%%%%%%%%%%%%%%%%%%%%%%%%%%%%%%%
%%%%%%%%%%%%%%%%%%%%%%%%%%%%%%%%%%%%%%%%%%%%%%%%%%%%%%%%%%%%%%%%%%%%%%%%%%%%%%%%%%%%%%%%%%%%%%%%%%%%%
%%%%%%%%%%%%%%%%%%%%%%%%%%%%%%%%%%%%%%%%%%%%%%%%%%%%%%%%%%%%%%%%%%%%%%%%%%%%%%%%%%%%%%%%%%%%%%%%%%%%%

%%%%%%%%%%%%%%%%%%%%%%%%%%%%%%%%%%%%%%%%%%%%%%%%%%%%%%%%%%%%%%%%%%%%%%%%%%%%%%%%%%%%%%%%%%%%%%%%%%%%%
%%%%%%%%%%%%%%%%%%%%%%%%%%%%%%%%%%%%%%%%%%%%%%%%%%%%%%%%%%%%%%%%%%%%%%%%%%%%%%%%%%%%%%%%%%%%%%%%%%%%%
%%%%%%%%%%%%%%%%%%%%%%%%%%%%%%%%%%%%%%%%%%%%%%%%%%%%%%%%%%%%%%%%%%%%%%%%%%%%%%%%%%%%%%%%%%%%%%%%%%%%%

\head\sectProofMainTh. Uniqueness of $MW$-curves up to isotopy (proof of Theorem 1)
\endhead

Let $K$ be as in Theorem 1. So, $K$ is a smooth rational curve in $\P^3$ of degree $d\ge 3$ and $w(K)=N_d$.

Given a point $p\in\P^3$, let $\pi_p:\P^3\setminus\{p\}\to\P^2$ be the projection from $p$ and
%(here $\P^2$ can be considered as the projectivization of the tangent space $T_p(\P^3)$).
let $\hat\pi_p:K\to\P^2$ be the restriction of $\pi_p$ to $K$. 
If $p\in K$, then we extend $\hat\pi_p$ to $p$
by continuity, thus $\pi_p^{-1}(\hat\pi_p(p))=T_p$
where $T_p$ is the tangent line to $K$ at $p$.

If $p\in\RP^3$, then we set $C_p=\hat\pi_p(K)$. Note that
$$
   \deg C_p = \cases d, & p\not\in K\\ d-1, & p\in K.\endcases
$$

Recall that an algebraic curve $C$ (maybe, singular)
of degree $m$ in $\RP^2$ is called {\it hyperbolic} with respect to a point $q\in\RP^2$
(which may or may not belong to $C$),
if any real line through $q$ intersects $C$ at $m$ real points counting the multiplicities.
We denote:
$$
   \hyp(C) = \{q\mid\text{$C$ is hyperbolic with respect to $q$}\}.            \eqno(\eqDefHyp)
$$
It is easy to check that $\hyp(C)$ is either empty or a convex closed set. It is possible that $\hyp(C)$
contains only one point. In this case, the point should be singular. For example,
if $C$ is a cuspidal cubic, then $\hyp(C)$ consists of the cusp only.

Similarly, we say that $K$ is hyperbolic with respect to a real line $L$ if, for
any real plane $P$ passing through $L$, each intersection point of $K$ and $P\setminus L$ is real.

The following two properties of hyperbolic curves are immediate from the definition:

\proclaim{ Lemma \lemHypC } Let $C$ be a real plane curve, $q\in\hyp(C)$, and
$q_1\in\R C\setminus\{q\}$.
Then each local branch of $C$ at $q_1$ is
smooth, real, and transverse to the line $(qq_1)$. The projection from $q$
defines a covering $\R\tilde C\to\RP^1$ where $\tilde C$ is the normalization of~$C$.
\qed
\endproclaim

\proclaim{ Lemma \lemHypKC }
Let $L$ be a real line in $\P^3$ and $p\in\R L$. Then
$K$ is hyperbolic with respect to $L$ if and only if
$\pi_p(L)\in\hyp(C_p)$. \qed
\endproclaim

%Using the arguments as in [\refMO; \S5], one easily obtains the following

\proclaim{ Lemma \lemSmooth } If $p\in\R K$, then $\hat\pi_p(p)$ is a smooth
point of $\R C_p$ and $T_p\cap K=\{p\}$.
\endproclaim

\demo{ Proof } By Proposition 1 %(see Section \sectTors),
the torsion of $\R K$
does not vanish at $p$, hence
the image of the germ $(K,p)$ under the projection $\hat\pi_p$ is
a smooth local branch of $C_p$ at $q=\hat\pi_p(p)$. Suppose that $C_p$ has another local
branch at $q$ which is the projection of a germ $(K,p_1)$.
Let $p_0$ be a real point on the line $(pp_1)$, $p_0\not\in\{p,p_1\}$.
Then $C_{p_0}$ has (at least) two branches at $\pi_{p_0}(p)=\pi_{p_0}(p_1)$ one of whom
is cuspidal. In this case by perturbing $K$ we may obtain either a non-real crossings or
a pair of real crossings of opposite signs which
contradicts the maximality of $w(K)$ (cp. the end of the proof of Proposition 1).
\qed\enddemo

\proclaim{ Lemma \lemHyperb }
If $p\in\R K$, then $\hyp(C_p)$ is the closure of a component of $\RP^2\setminus\R C_p$
and $\hat\pi_p(p)$ is a smooth point of its boundary.
\endproclaim

\demo{ Proof }
Let  $q=\hat\pi_p(p)$. It is a
smooth point of $C_p$ by Lemma \lemSmooth.
For a real singular point $u$ of $C_p$, let $\sigma(u)$ be its contribution to $w(K)$,
i.e. the sum of the signs of crossings of a nodal perturbation of $u$.
%We assume for simplicity that
%all singularities of $C_p$ are ordinary double points
%To simplify the notation, we assume that
%all singularities of $C_p$ are ordinary double points
%(see Remark 1 below).
Then, similarly to [\refMO; Proposition 21], we have
$$
   w(K) = i(q') + i(q'') + \sum_{u\in\Sing(\R C_p)} \sigma(u),            \eqno(\eqIndex)
$$
where $q',q''\not\in C_p$ are points close to $q$ on different sides of $\R C_p$, and
$i(x)$ for $x\not\in C_p$ is one half of the image of $[\R C_p]$ under
the isomorphism $H_1(\RP^2\setminus\{x\})\cong\Z$
(see [\refMO; \S6] for the choice of the orientations).
%and $\sigma(u)=\pm1$ are the same signs as in the definition of $w$.
It is clear that $i(q')+i(q'')\le\deg C_p-1=d-2$ and
$\sum_u\sigma(u)\le\card\Sing(\R C_p)\le (d-2)(d-3)/2$. The sum of these two upper bounds is $N_d$,
thus $w(K)=N_d$ implies the equality sign in the both estimates.
It remains to note that $i(q') + i(q'')=d-2$ implies $q\in\hyp(C_p)$.

The fact that $\hyp(C_p)$ is the closure of a component of the complement of $C_p$ follows
from the discussion after (\eqDefHyp) because $q$ is smooth on $C_p$. \qed
\enddemo

Recall that the tangent line to $K$ at $p\in K$ is denoted by $T_p$.
Let us set
$$
       T = \bigcup_{p\in\R K} \R T_p
$$

\proclaim{ Lemma \lemHypK }
Suppose that $K$ is hyperbolic with respect to a real line $L$ and let
$p\in(\R K)\setminus L$. Then $L\cap T_p=\varnothing$.
\endproclaim

\demo{ Proof } 
Combine Lemma \lemHypC\ and Lemma \lemHypKC.
\qed\enddemo

\proclaim{ Lemma \lemDisjoint } Let $p_1$ and $p_2$ be two distinct points on $\R K$.
Then $T_{p_1}\cap T_{p_2}=\varnothing$. 
\endproclaim

\demo{ Proof }
Let $L=T_{p_1}$. Then
$K$ is hyperbolic with respect to $L$ by Lemma \lemHyperb\ combined with
Lemma \lemHypKC, and we have $p_2\not\in L$ by Lemma \lemSmooth.
Hence the result follows from Lemma \lemHypK.
\qed\enddemo

Thus $T$ is a disjoint union of a continuous family of real projective lines (topologically, circles)
parametrized by $\R K$.
We are going to show that the pair $(T,\R K)$  is isotopic in $\RP^3$
to a hyperboloid with a projective torus link
$T_{proj}(d,d-2)$ sitting in it.
Note that $T$ is not smooth. It has a cuspidal edge along $\R K$. %(cp.~Remark \remCycloid).

\proclaim{ Lemma \lemTwoLines } There exist two real lines $L_1$ and $L_2$ such that $K$ is
hyperbolic with respect to each of them, $L_1\cap K=\varnothing$, and
$L_2$ crosses $K$ without tangency at a pair of complex conjugate points.
\endproclaim

\midinsert
\epsfxsize=60mm
\centerline{\epsfbox{ 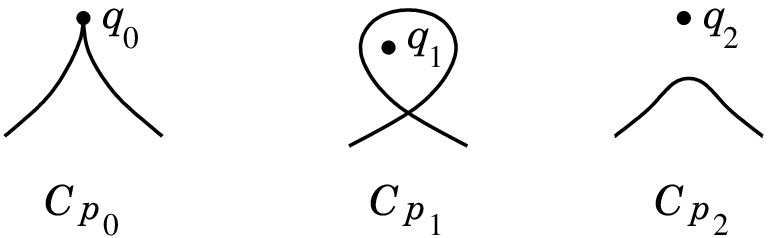 }}
\botcaption{ Figure \figCusp } Two perturbations of a cusp in the proof of Lemma \lemTwoLines
\endcaption
\endinsert

\demo{ Proof } 
Let us choose a point $p\in\R K$ % such that $T_p\cap K=\{p\}$ 
and let $p_0\in\R T_p\setminus\{p\}$.
Then $\pi_p(T_p)=\hat\pi_p(p)\in\hyp(C_p)$ by Lemma \lemHyperb\ whence $K$ is hyperbolic with
respect to $T_p$ by Lemma \lemHypKC. Let $q_0=\pi_{p_0}(p)$. Then, again by Lemma \lemHypKC, we have
$q_0\in\hyp(C_{p_0})$. The curve $C_{p_0}$ has a cusp at $q_0$ because the torsion at $p$ is nonzero.
Let $p_1$ and $p_2$ be points close to $p_0$ and chosen on different sides of the
osculating plane of $\R K$ at $p$. Then $C_{p_1}$ and $C_{p_2}$ are
obtained from $C_{p_0}$ by a perturbation of the cusp as shown in Figure \figCusp\ where
$q_2$ is a solitary node of $C_{p_2}$ (a point where two complex conjugate branches cross).
Then we set $L_j = \pi_{p_j}^{-1}(q_j)$, $j=1,2$, where the points $q_1$ and $q_2$ are chosen
as in Figure \figCusp. The fact that $q_0\in\hyp(C_{p_0})$ implies $q_j\in\hyp(C_{p_j})$, $j=1,2$,
whence the hyperbolicity of $K$ with respect to $L_j$ by Lemma \lemHypKC.
\qed\enddemo

\demo{ Proof of Theorem 1 }
Let $L_1$ and $L_2$ be as in Lemma \lemTwoLines.

The line $L_1$ is disjoint from $T$ by Lemma \lemHypK.
Let $P$ be a real plane through $L_1$. 
Again by Lemma \lemHypK, $P$ crosses each line $T_p$, $p\in\R K$, at a single point.
Let us denote this point by $\xi_P(p)$. Then $\xi_P:\R K\to\R P$ is a continuous mapping.
It is injective by Lemma \lemDisjoint\ and its image (which is $T\cap\R P$) is disjoint from $L_1$.
Hence $T\cap\R P$ is a Jordan curve in the affine real plane $\R P\setminus L_1$.
Let $D_P$ be the disk bounded by this Jordan curve and let $U_1=\bigcup_P D_P$ where
$P$ runs through all the real planes through $L_1$.
Then $U_1$ is fibered by disks over a circle which parametrizes the pencil of planes through $L_1$.
Since $\RP^3$ is orientable, this fibration is trivial, thus $U_1$ is a solid torus and $\partial U_1=T$.
Each $P$ transversally crosses $K$ at $d$ real points, thus $\R K$ sits in $T$ and
it realizes the homology class $d\alpha$ where $\alpha$ is a generator of $H_1(U_1)$.

The same arguments applied to the line $L_2$ show that $T$ bounds a solid torus $U_2$
such that  $\R K$ realizes the homology class $(d-2)\beta$ where $\beta$ is a generator of $H_1(U_2)$.
We conclude that the lift of $K$ on $\sph^3$ is $T(d,d-2)$ and the result follows.
\qed\enddemo

We see that $T$ cuts $\RP^3$ into two solid tori $U_1$ and $U_2$
such that $\R L_1\subset U_2$ and $\R L_2\subset U_1$.

\proclaim{ Proposition \propU } {\rm(Compare with Corollary 1).}
Let $p$ be a generic point of $\RP^3$. Then $C_p$ has only real double points.
If $p\in U_1$, then all the double points have real local branches and the interior of $\hyp(C_p)$ is non-empty.
If $p\in U_2$, then one double point $q$ is solitary (i.e. has complex conjugate local branches),
all the other double points have real local branches, and $\hyp(C_p)=\{q\}$.
\endproclaim

\demo{ Proof } Let us consider a generic path $p(t)$ which relate the given point with a point on $T$.
It defines a continuous deformation of the knot diagram which is a sequence of Reidemeister moves
(R1) -- (R3). However, (R2) is impossible because it involves a negative crossing and (R1) smay
occur only when $p(t)$ passes through $T$. Thus the number and the nature of double points
does not change during the deformation. The projection from a point of $T$ is cuspidal and it is
hyperbolic with respect to the cusp, so the result follows from Lemma \lemHypC.

%The assertion $\hyp(C_p)\ne\varnothing$ 
Non-emptiness of the interior of $\hyp(C_p)$
in case $p\in U_1$, follows from the fact that $\hyp(C_p)$
can disappear only by a move (R3). This is however impossible because all crossings are positive
and the boundary orientation on $\partial(\hyp(C_p))$ agrees with an orientation of $\R C_p$
due to Lemma \lemHypC\ (see Figure \figR).
\qed\enddemo

\midinsert
\centerline{\lower-10mm\hbox{$\hyp(C_p)$}\epsfxsize=20mm\epsfbox{ 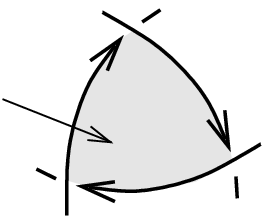 }}
\botcaption{ Figure \figR } Impossibility of a move (R3) which eliminates $\hyp(C_p)$
\endcaption
\endinsert

%%%%%%%%%%%%%%%%%%%%%%%%%%%%%%%%%%%%%%%%%%%%%%%%%%%%%%%%%%%%%%%%%%%%%%%%%%%%%%%%%%%%%%%%%%%%%%%%%%%%%
%%%%%%%%%%%%%%%%%%%%%%%%%%%%%%%%%%%%%%%%%%%%%%%%%%%%%%%%%%%%%%%%%%%%%%%%%%%%%%%%%%%%%%%%%%%%%%%%%%%%%
%%%%%%%%%%%%%%%%%%%%%%%%%%%%%%%%%%%%%%%%%%%%%%%%%%%%%%%%%%%%%%%%%%%%%%%%%%%%%%%%%%%%%%%%%%%%%%%%%%%%%

%%%%%%%%%%%%%%%%%%%%%%%%%%%%%%%%%%%%%%%%%%%%%%%%%%%%%%%%%%%%%%%%%%%%%%%%%%%%%%%%%%%%%%%%%%%%%%%%%%%%%
%%%%%%%%%%%%%%%%%%%%%%%%%%%%%%%%%%%%%%%%%%%%%%%%%%%%%%%%%%%%%%%%%%%%%%%%%%%%%%%%%%%%%%%%%%%%%%%%%%%%%
%%%%%%%%%%%%%%%%%%%%%%%%%%%%%%%%%%%%%%%%%%%%%%%%%%%%%%%%%%%%%%%%%%%%%%%%%%%%%%%%%%%%%%%%%%%%%%%%%%%%%

\subhead Acknowledgement
\endsubhead
We are grateful to Oleg Viro for very useful discussions.

\enddocument